\documentclass{amsproc}
\usepackage{enumerate}
\usepackage{graphicx}
\usepackage{amssymb}
\usepackage{epstopdf}
\usepackage{amsmath,amsthm,amscd,amssymb}
\usepackage{latexsym}
\usepackage[colorlinks,citecolor=red,pagebackref,hypertexnames=false]{hyperref}
\usepackage{geometry}                
\numberwithin{equation}{section}

\theoremstyle{plain}
\newtheorem{theorem}{Theorem}[section]
\newtheorem{lemma}[theorem]{Lemma}
\newtheorem{corollary}[theorem]{Corollary}

\newtheorem{conjecture}[theorem]{Conjecture}

\theoremstyle{definition}
\newtheorem{definition}[theorem]{Definition}

\newtheorem{case[theorem]}{Case}

\theoremstyle{remark}

\numberwithin{equation}{section}

\def\dH{\dim_{{\mathcal H}}}
\def\R{\Bbb R}

\def\e{\epsilon}

\begin{document}

\title{Pinned distance problem, slicing measures and local smoothing estimates} 


\author{Alex Iosevich and Bochen Liu}

\date{today}

\email{iosevich@math.rochester.edu}
\address{Department of Mathematics, University of Rochester, Rochester, NY}

\email{bochen.liu@rochester.edu}
\address{Department of Mathematics, University of Rochester, Rochester, NY}

\thanks{This work was partially supported by the NSA Grant H98230-15-1-0319}

\maketitle
\begin{abstract} We improve the Peres-Schlag result on pinned distances in sets of a given Hausdorff dimension. In particular, for Euclidean distances, with $$\Delta^y(E) = \{|x-y|:x\in E\},$$ we prove that for any $E, F\subset\R^d$, there exists a probability measure $\mu_F$ on $F$ such that for $\mu_F$-a.e. $y\in F$, 

\vskip.125in

\begin{itemize}
   \item $\dH(\Delta^y(E))\geq\beta$ if $\dH(E)+\frac{d-1}{d+1}\dH(F)>d-1+\beta$;
   \item $\Delta^y(E)$ has positive Lebesgue measure if $\dH(E)+\frac{d-1}{d+1}\dH(F)>d$;
   \item $\Delta^y(E)$ has non-empty interior if $\dH(E)+\frac{d-1}{d+1}\dH(F)>d+1$.
\end{itemize}

\vskip.125in

We also show that in the case when $\dH(E)+\frac{d-1}{d+1}\dH(F)>d$, for $\mu_F$-a.e. $y\in F$,
$$ \left\{t\in\R : \dH(\{x\in E:|x-y|=t\}) \geq \dH(E)+\frac{d+1}{d-1}\dH(F)-d \right\}   $$
has positive Lebesgue measure. This describes dimensions of slicing subsets of $E$, sliced by spheres centered at $y$.

In our proof, local smoothing estimates of Fourier integral operators (FIO) plays a crucial role. In turn, we obtain results on sharpness of local smoothing estimates by constructing geometric counterexamples. 
\end{abstract}

\section{Introduction}

\subsection{Distance problem and pinned distance problem}
One of the most important open problems in geometric measure theory is the Falconer distance conjecture (\cite{Fal85}). It says that given any set $E\subset\R^d$, $d\geq 2$, its distance set
$$\Delta(E)=\{|x-y|: x,y\in E\}$$
has positive Lebesgue measure whenever $\dH(E)>\frac{d}{2}$. The best currently known results are due to Wolff in $\R^2$ (\cite{Wol99}) and Erdogan in higher dimensions (\cite{Erd05}). They proved that $\Delta(E)$
has positive Lebesgue measure whenever $\dH(E)>\frac{d}{2}+\frac{1}{3}$. There are also results where $E$ is assumed to have special structures (\cite{IL16}, \cite{Orp17}, \cite{Shm16}).

\vskip.125in

More generally, we may ask how large the Hausdorff dimensions of $E\subset\R^d$ need to be to ensure the set
$$\Delta_\Phi(E)=\{\Phi(x,y): x,y\in E\}$$
has positive Lebesgue for suitable functions $\Phi : \R^d\times\R^d\rightarrow\R$. In \cite{EIT11}, Eswarathasan, Iosevich and Taylor prove that if $\Phi:\R^d\times\R^d\rightarrow\R$ has non-zero Monge-Ampere determinant (also called Phong-Stein rotational curvature condition), i.e.
\begin{equation}\label{PS}
det
\begin{pmatrix}
  0&\nabla_x\Phi\\\nabla_y\Phi&\frac{\partial^2\Phi}{\partial x\partial y}
\end{pmatrix}
\neq 0,
\end{equation}
then $\Delta_\Phi(E)$ has positive Lebesgue measure whenever $\dH(E)>\frac{d+1}{2}$. 

A corollary of this result is that if $E$ is a closed subset of a closed compact Riemannian manifold $M$ of dimension $d \ge 2$ of Hausdorff dimension $>\frac{d+1}{2}$, then $\Delta_{\rho}(E)$ has positive Lebesgue measure, where $\rho$ denotes the Riemannian metric on $M$. 

\vskip.125in

Another interesting version of the Falconer distance problem is the pinned distance problem. Given $E\subset\R^d$, we ask whether for "many" points $y\in \R^d$, the pinned distance set 
$$\Delta^y(E) = \{|x-y|:x\in E\}$$ 
has positive Lebesgue measure. This problem was first studied by Peres and Schlag (\cite{PS00}). They considered a very large class of functions $\Phi(x,y)$ that generalize orthogonal projections $\Phi(x,e)=x\cdot e$, $ e\in S^{d-1}$. For Euclidean distances, one can take $\Phi(x,y)=|x-y|$ and their result implies the following.

\begin{theorem}[Peres, Schlag, (2000) (\cite{PS00})]
\begin{equation}
   \label{PS1}
   \dH( \{y\in\R^d: \dH(\Delta^y(E))<\beta \})\leq d+\beta -\max\{\dH(E),1\},
\end{equation} 
\begin{equation}
   \label{PS2}
   \dH( \{y\in\R^d: |\Delta^y(E)| = 0 \})\leq d+1-\dH(E),
\end{equation} 
\begin{equation}
   \label{PS3}
   \dH( \{y\in\R^d: \text{Int}(\Delta^y(E)) = \emptyset\})\leq d+2-\dH(E).
\end{equation} 
\end{theorem}

In 2016, Iosevich, Taylor and Uriarte-Tuero (\cite{ITU16}) gave a straightforward proof of \eqref{PS2} for functions $\Phi(x,y)$ satisfying the Phong-Stein rotational curvature condition \eqref{PS} in the range $\dH(E)>\frac{d+1}{2}$. Shmerkin (\cite{shm17}) recently proved that
$$\dH( \{y\in\R^2: \dH(\Delta^y(E))<1\})\leq 1$$
for any $E\subset\R^2$ such that $\dH(E)=\dim_P(E)>1$, where $\dim_{P}(E)$ denotes the packing dimension of $E$. 

\vskip.125in

It is interesting to note that the Peres-Schlag bound is in general sharp. One can simply take $\Phi(x,y)=x\cdot y$ and $E=A\times\R^{k}\times\{0\}^{d-k-1}$ where $A\subset\R$ has Hausdorff dimension $\beta$. Then $\dH(E)=\beta+k$ while
$$\{y\in\R^d: \dH(\Delta^y(E))<\beta \}\supset\R\times\{0\}^{k}\times\R^{d-k-1},$$
whose Hausdorff dimension is $d-k=d+\beta-\dH(E)$. Therefore, if we want to improve Peres-Schlag's bound for, say Euclidean distances, we need a stronger assumption that rules out the case $\Phi(x,y)=x\cdot y$.

\vskip.125in

Notice that the key difference between $|x-y|$ and $x\cdot y$ is that if $y$ is fixed, $\{x:|x-y|=1\}$ is a sphere with non-zero Gaussian curvature, while $\{x:x\cdot y =1\}$ is just a hyperplane. With this in mind, we obtain the following improvement of the aforementioned Peres-Schlag's result when $\dH(E)\geq\frac{d+1}{2}$.

\begin{definition}
  We say that $\Phi\in C^\infty(\R^d\times\R^d)$ satisfies Sogge's cinematic curvature condition (\cite{Sog91}) if 
\begin{equation}\label{ccc}
  \text{for any $t>0$, $x\in \R^d$, $\{\nabla_y\Phi : \Phi(x,y)=t\}$
has nonzero Gaussian curvature.}\end{equation}
\end{definition}
\vskip.125in
\begin{theorem}
  \label{main}
Given $\Phi\in C^\infty(\R^d\times\R^d)$, $E, F\subset\R^d$. Suppose $\Phi$ satisfies the Phong-Stein rotational curvature condition \eqref{PS} and the cinematic curvature condition \eqref{ccc}. Then there exists a probability measure $\mu_F$ on $F$ such that for $\mu_F$ a.e. $y\in F$,
\begin{equation}
   \dH(\Delta_{\Phi}^y(E))\geq\beta \ \text{if}\ \dH(E)+\frac{d-1}{d+1}\dH(F)>d-1+\beta;
\end{equation} 
\begin{equation}
   |\Delta_{\Phi}^y(E)|>0\ \text{if}\ \dH(E)+\frac{d-1}{d+1}\dH(F)>d;
\end{equation} 
\begin{equation}
   \text{Int}(\Delta_{\Phi}^y(E)) \neq \emptyset\ \text{if}\ \dH(E)+\frac{d-1}{d+1}\dH(F)>d+1.
\end{equation} 
\end{theorem}
\vskip.125in
Improvement on pinned distance problem then follows easily.
\begin{corollary}
  \label{pinD}
Given $\Phi\in C^\infty(\R^d\times\R^d)$, $E\subset\R^d$. Suppose $\Omega\subset\R^d$ and $\Phi$ satisfies the Phong-Stein rotational curvature condition \eqref{PS} and the cinematic curvature condition \eqref{ccc} on $E\times\Omega$. Then
\begin{equation}\label{pinD1}
   \dH( \{y\in\Omega: \dH(\Delta_{\Phi}^y(E))<\beta \})\leq d+1+\frac{d+1}{d-1}\beta -\frac{d+1}{d-1}\dH(E),
\end{equation} 
\begin{equation}\label{pinD2}
   \dH( \{y\in\Omega: |\Delta_{\Phi}^y(E)| = 0 \})\leq \frac{d(d+1)}{d-1} -\frac{d+1}{d-1}\dH(E),
\end{equation} 
\begin{equation}\label{pinD3}
   \dH( \{y\in\Omega: \text{Int}(\Delta_{\Phi}^y(E)) = \emptyset\})\leq \frac{(d+1)^2}{d-1} -\frac{d+1}{d-1}\dH(E).
\end{equation} 
\end{corollary}

One can check that \eqref{pinD1} improves \eqref{PS1} when $\dH(E)>\frac{d-1}{2}+\beta$; \eqref{pinD2} improves \eqref{PS2} when $\dH(E)>\frac{d+1}{2}$ and \eqref{pinD3} improves \eqref{PS3} when $\dH(E)>\frac{d+3}{2}$.

\vskip.125in 

In particular, when $\dH(E)=d$, we have the following sharp corollary.
\begin{corollary}
  Suppose $E\subset\R^d$, $\dH(E)=d$. With notations above,
$$\dH( \{y\in\Omega: |\Delta_{\Phi}^y(E)| = 0 \})=0.$$
\end{corollary}
\vskip.25in
\subsection{Dimensions of slicing sets}
Given a family of functions $\{\pi_\lambda\}_{\lambda\in \Lambda}$ and $E\subset\R^d$ such that $\pi_\lambda(x)\neq\pi_{\lambda'}(x)$ whenever $\lambda\neq\lambda'$, then $\{\pi^{-1}_\lambda(t)\}_{t\in\R}$ slices $E$ into pieces and we can study not only the size of images $\{\pi_\lambda(x):\lambda\in\Lambda, x\in E\}$, but also the size of those slices $E\cap\{\pi^{-1}_\lambda(t)\}$. The most classical result is due to Marstrand.

\begin{theorem}[Marstrand, 1954] Suppose $E\subset\R^2$ is a Borel set and let $\pi_\theta(x)=x\cdot (\cos\theta,\sin\theta)$. Then for almost all $\theta\in [0,2\pi)$,
\begin{equation}\label{proj}
\begin{aligned}
\dH(\pi_\theta(E))=\dH(E) \ \text{if}\  \dH(E)\leq 1;\\
|\pi_\theta(E)|>0 \ \text{if}\  \dH(E)>1.
\end{aligned}
\end{equation}
Moreover, when $\dH(E)>1$, the typical lines with direction $\theta$ which intersect $E$ intersect it in dimension $\dH(E)-1$, i.e.,
\begin{equation}\label{slice}
|\{t\in\R:\dH(\pi^{-1}_\theta(t)\cap E)=\dH(E)-1\}|>0.
\end{equation}
\end{theorem}

A variety of generalizations of this result have been obtained over the years. See, for example, \cite{Mat14} and the references contained therein. 

\vskip.125in

What Peres and Schlag proved in \cite{PS00} generalized \eqref{proj} from orthogonal projections to a large class of functions so-called generalized projections, including, in particular, Euclidean distances $\pi_\lambda (x)= |x-\lambda|$, $x,\lambda\in\R^d$. In \cite{Orp14}, Orponen generalized \eqref{slice} from orthogonal projections to generalized projections. 

\vskip.125in

Since Theorem \ref{main} improves Peres-Schlag's result on dimensions of images of projections, it is natural to ask whether an improvement on dimensions of slicing sets can also be obtained. Its geometric meaning is that under the assumptions of Theorem \ref{main}, a typical distance $t$ occurs statistically often. 

\begin{theorem}\label{slicingthm}
Under the assumptions of Theorem \ref{main}, with $E, F\subset\R^d$ and  
$$\dH(E)+\frac{d+1}{d-1}\dH(F)>d,$$ there exists a probability measure $\mu_F$ on $F$ such that for $\mu_F$-a.e. $y\in F$,
$$|\{\Phi(x,y):x\in E\}|>0  $$
and 
$$|\{t : \dH(\{x\in E:\Phi(x,y)=t\}) \geq \dH(E)+\frac{d+1}{d-1}\dH(F)-d \}|>0.     $$

\end{theorem}

\vskip.125in

\subsection{Local smoothing of Fourier integral operators}
The study of Fourier integral operators arises in the study of the wave equation,
\begin{equation*}
  \begin{cases}
    \frac{\partial^2 u}{\partial t^2}=\Delta u,\\
    u|_{t=0}=f,\  \frac{\partial u}{\partial t}|_{t=0}=0,
  \end{cases}
\end{equation*}
whose solution is given by the real part of 
$$\mathcal{F} f (y,t):= \int e^{2\pi i y\cdot\xi} e^{-2\pi i t|\xi|}\hat{f}(\xi)\,d\xi. $$

More generally, we may consider the Fourier integral operators (FIO), introduced by H\"ormander in 1970s (\cite{Hor71}). In this paper, we only present a simple form of FIO for convenience. For general information, see e.g. \cite{Sog93} and the references therein.
\begin{definition}\label{FIO}
Suppose $\Phi\in C^\infty(U\times V)$, where $U, V\subset\R^d$ are open domains. For any $f\in C_0^\infty(V)$, one can define
\begin{equation*}
\begin{aligned}
\mathcal{F} f(y,t) =& \int_{\R^d} \left(\int_{\R} e^{-2\pi i (\Phi(x,y)-t)\tau}\,d\tau\right)f(x)\,dx\\=&\lim_{N \rightarrow \infty}\int_{\R^d} \left(\int^N_{-N} e^{-2\pi i (\Phi(x,y)-t)\tau}\,d\tau\right)f(x)\,dx.
\end{aligned}
\end{equation*} 
\end{definition}

\vskip.125in

Throughout this paper, we use the notation 
$$\alpha_{d,p,\gamma}=-\frac{d-1}{2}+\gamma+(d-1)\left|\frac{1}{2}-\frac{1}{p}\right|.$$

\vskip.125in

It is known that if $\Phi$ satisfies the Phong-Stein rotational curvature condition \eqref{PS}, then for any fixed $t>0$, 
$$ (I-\Delta)^{\frac{\gamma}{2}}\mathcal{F}(\cdot,t):L^p_{\alpha_{d,p,\gamma}}(U) \rightarrow L^p(V),$$
for any $1<p<\infty$. This result is sharp if $t$ is fixed. But if $\Phi$ satisfies the cinematic curvature condition \eqref{ccc}, there is a gain of regularity for $2<p<\infty$ by taking average in $t\approx 1$. This phenomenon is called local smoothing.

\begin{theorem}[$L^2\rightarrow L^p$ local smoothing estimate, Seeger, Sogge, Stein, 1991]\label{l2} Let $\mathcal{F}$ be as in Definition \ref{FIO}. Suppose $\Phi$ satisfies the Phong-Stein rotational curvature condition \eqref{PS} and the cinematic curvature condition \eqref{ccc} on $U\times V$. Then 
\begin{equation}\label{l2lp1}\left|\left|\left(\int_1^2|(I-\Delta)^{\frac{\gamma}{2}}\mathcal{F} f(\cdot,t)|^2\,dt\right)^\frac{1}{2}\right|\right|_{L^p(V)}\lesssim ||f||_{L^2_{\alpha_{d,p,\gamma}-\delta_2(d,p)}(U)},
\end{equation}
where \begin{equation}\label{delta_2}
  \delta_2(d,p)=
  \begin{cases}
  \frac{1}{p}, & p\geq\frac{2(d+1)}{d-1};\\\frac{1}{2}(d-1)(\frac{1}{2}-\frac{1}{p}), & 2<p\leq\frac{2(d+1)}{d-1}.
  \end{cases} 
\end{equation}
\end{theorem}

It is known that \eqref{l2lp1} does not hold with $\delta_2(d,p)=\frac{1}{p}$ for any $p<\frac{2(d+1)}{(d-1)}$. Also $L^2_{\alpha_{d,p,\gamma}-\frac{1}{p}-\epsilon}(U)\rightarrow L^p(L^2([1,2]))(V)$ breaks down for any $\epsilon>0, p>2$.

\vskip.125in

It is more interesting to consider $L^p\rightarrow L^p$ local smoothing estimate because it is related to the Bochner-Riesz conjecture. Sogge's local smoothing conjecture (\cite{Sog91}) takes the following form in our setup. 

\begin{conjecture}[$L^p\rightarrow L^p$ local smoothing estimate, Sogge, 1991]\label{Sog}
  Let $\mathcal{F}$ be as in Definition \ref{FIO}. Suppose $\Phi$ satisfies the Phong-Stein rotational curvature condition \eqref{PS} and the cinematic curvature condition \eqref{ccc} on $U\times V$. Then 
  \begin{equation}\label{ls}
  (I-\Delta)^{\frac{\gamma}{2}}\mathcal{F}: L^p_{\alpha_{d,p,\gamma}-\delta_p(d,p)+\epsilon} (U)\rightarrow L^p_{}(V\times [1,2]) 
  \end{equation}
  for any $\epsilon>0$ with
\begin{equation}
  \delta_p(d,p)=
  \begin{cases}
  \frac{1}{p}, & p\geq\frac{2d}{d-1};\\(d-1)(\frac{1}{2}-\frac{1}{p}), & 2<p\leq\frac{2d}{d-1}.
  \end{cases}
\end{equation}
\end{conjecture}
The best currently known result, accumulating efforts by Sogge (\cite{Sog91}), Mockenhaupt, Seeger, Sogge (\cite{MSS93}) and Wolff (\cite{Wol00}), is due to Bourgain and Demeter (\cite{BD15}). They proved that \eqref{ls} holds with  
\begin{equation}
  \delta_p(d,p)=
  \begin{cases}
  \frac{1}{p}, & p\geq\frac{2(d+1)}{d-1};\\\frac{1}{2}(d-1)(\frac{1}{2}-\frac{1}{p}), & 2<p\leq\frac{2(d+1)}{d-1}. 
  \end{cases}
\end{equation}

It is known that $\delta_p(d,p) = \frac{1}{p}$ is the best possible for $p\geq\frac{2d}{d-1}$ due to the sharpness of Bochner-Riesz conjecture. There are also examples showing $\delta_p(d,p)>\frac{1}{p}$ fails for all $2<p<\infty$. In \cite{MS97}, Minicozzi and Sogge construct a metric $\Phi$ on odd dimensional spaces such that $\delta_p(d,p)=\frac{1}{p}$ works only if $p\geq\frac{2(d+2)}{d}$. Whether Sogge's conjecture is optimal in the range $p \ge \frac{2d}{d-1}$ is an open question. 

\vskip.125in

We will see that the proof of Theorem \ref{main} relies on local smoothing estimates. Therefore, in turn, geometric counterexamples of Theorem \ref{main} imply sharpness of local smoothing estimates. 
\begin{theorem}[Sharpness of Theorem \ref{l2}]\label{sharpl2}Suppose $\Phi$ satisfies assumptions in Theorem \ref{l2}.
  \begin{enumerate}[(i)]
    \item For all $\Phi$, 
    $$L^2_{\alpha_{d,p,\gamma}-\frac{1}{p}-\epsilon}(U)\rightarrow L^p(L^2([1,2]))(V)$$ breaks down for any $\epsilon>0, p>2$.
    \item In odd dimensions, there exists $\Phi$ to show that Theorem \ref{l2} is sharp for all $p>2$.
    \item In even dimensions, there exists $\Phi$ such that \eqref{l2lp1} breaks down with any $\delta_2(d,p)>\frac{1}{2}d(\frac{1}{2}-\frac{1}{p})$
  \end{enumerate}
\end{theorem}
\vskip.125in
\begin{theorem}[Sharpness of Conjecture \ref{Sog}]\label{sharplp}Suppose $\Phi$ satisfies assumptions in Conjecture \ref{Sog}
  
  \begin{enumerate}[(i)]
    \item For all $\Phi$, 
    $$L^p_{\alpha_{d,p,\gamma}-\frac{1}{p}-\epsilon}(U)\rightarrow L^p(V\times [1,2])$$ breaks down for any $\epsilon>0, p>2$.
    \item In odd dimensions, there exists $\Phi$ to show that Conjecture \ref{Sog} is the best possible for all $2<p\leq\frac{2d}{d-1}$.
    \item In even dimensions, there exists $\Phi$ such that \eqref{ls} breaks down with any $\delta_p(d,p)> d(\frac{1}{2}-\frac{1}{p})$.
  \end{enumerate}
\end{theorem}

\subsection*{Notation}
\begin{itemize}
  \item $s_E$ denotes the Hausdorff dimension of $E\subset\R^d$.
  \item $X\lesssim Y$ means there exists a constant $C$ such that $X\leq C Y$. $X\lesssim_\epsilon Y$ means the implicit constant only depends on $\epsilon$.
  \item $||\nu||^2_{L^2_\gamma}:=\int|\hat{\nu}(\xi)|^2|\xi|^{2\gamma}\,d\xi.$
  \item $\alpha_{d,p,\gamma}=-\frac{d-1}{2}+\gamma+(d-1)\left|\frac{1}{2}-\frac{1}{p}\right|$.
\end{itemize}
\vskip.125in
{\bf Acknowledgements.} The second listed author would like to thank Professor Ka-Sing Lau for the financial support of research assistantship in Chinese University of Hong Kong.

\vskip.125in

\section{Proof of Theorem \ref{main}}
For any $E\subset\mathbb{R}^n$, there exists Frostman measures supported on it that reflects its Hausdorff dimension. 
\begin{lemma}[Frostman Lemma, see e.g. \cite{Mat95}]\label{Frostmanlemma}
Denote $\mathcal{H}^s$ as the $s$-dimensional Hausdorff measure and $E\subset\mathbb{R}^n$. Then $\mathcal{H}^s(E)>0$ if and only if there exists a probability measure $\mu$ on $E$, such that 
$$\mu(B(x,r))\lesssim r^s, \ \forall\,x\in \mathbb{R}^n,\ \forall\,r>0.$$
\end{lemma}
\vskip.125in
The definition of Hausdorff dimension states that $\dH(E)=\sup\{s:\mathcal{H}^s(E)>0\}$, denoted by $s_E$. So by Lemma \ref{Frostmanlemma}, for any $E\subset\R^d$ and $\epsilon>0$, there exists $\mu_E$ on $E$ such that 
\begin{equation}\label{Frostman}\mu_E(B(x,r))\lesssim_\epsilon r^{s_E-\epsilon}, \ \forall\,x\in \mathbb{R}^n,\ \forall\,r>0.
\end{equation}
\vskip.125in
To study the size of the support of a measure $\nu$, we need the following lemma.
\begin{lemma}[see, e.g. Theorem 5.4 in \cite{Mat15}]\label{supp}
Suppose $\nu$ is a probability measure on $\R$ and $||\nu||_{L^2_{\gamma}}<\infty$.

\begin{enumerate}[(i)]
  \item If $\gamma<0$, the support of $\nu$ has Hausdorff dimension at least $1+2\gamma$.
  \item If $\gamma\geq 0$, the support of $\nu$ has positive Lebesgue measure.
  \item If $\gamma>\frac{1}{2}$, the support of $\nu$ has non-empty interior.
\end{enumerate}
\end{lemma}
\vskip.125in
Fix $y\in\R^d$, define a measure $\nu_y$ on
$$\Delta^y_{\Phi}(E)=\{\Phi(x,y):x\in E\}$$
by
\begin{equation}\label{pinmeasure}
\begin{aligned}
\int f(t)d\nu_y(t) = &\int f(\Phi(x,y))\,d\mu_E(x)\\
=& \iint \e^{2\pi i \Phi(x,y)\tau}\hat{f}(\tau)\,d\tau\,d\mu_E(x)\\
=& \iint \e^{2\pi i \Phi(x,y)\tau}\left(\int e^{-2\pi i t\tau}f(t)\,dt\right)\,d\tau\,d\mu_E(x).\\
\end{aligned}
\end{equation}
Therefore, in the sense of distribution,
$$d\nu_y(t)=\int \left(\int e^{-2\pi i (\Phi(x,y)-t)\tau}\,d\tau\right)\,d\mu_E(x),$$ 
which will be denoted as $\mathcal{F}\mu(y,t)$ as in Definition \ref{FIO}.
\vskip.125in

By Lemma \ref{supp}, to prove Theorem \ref{main}, it suffices to show that for any $F\subset\R^d$ of Hausdorff dimension $s_F>d+1+ (1+2\gamma)\frac{(d+1)}{d-1}-\frac{d+1}{d-1}s_E$ and for any Frostman measure $\mu_F$ on $F$ in \eqref{Frostman},
\begin{equation}\label{goal}
\int ||\nu_y||_{L^2_{\gamma}}\,d\mu_F(y)=\int ||\partial_t^\gamma\mathcal{F}\mu(y,t)||_{L^2(dt)}\,d\mu_F(y)<\infty.
\end{equation}

We first decompose $\mu_E,\mu_F$ into Littlewood-Paley pieces. There exists $\phi(\xi)\in C_0^\infty(\R^d)$, supported on $\{\frac{1}{2}<|\xi|<2\}$ such that

$$\sum_{j=1}^\infty \phi \left(\frac{\xi}{2^j}\right) = 1,\ \forall |\xi|\geq 1.   $$

Denote $\phi_j(\cdot)=\phi(\frac{\cdot}{2^j})$, $\phi_0=1-\sum_{j=1}^\infty {\phi_j}$, $\widehat{\mu_E^j}=\widehat{\mu_E}\phi_j$, $\widehat{\mu_F^j}=\widehat{\mu_F}\phi_j$, then
$$\mu_E = \sum^\infty_{j=0} \mu_E^j,\ \mu_F = \sum^\infty_{j=0} \mu_F^j,$$
and
\begin{equation}\label{decomp}
\begin{aligned}
\int ||\partial_t^\gamma\mathcal{F}\mu_E(y,t)||_{L^2(dt)}\,d\mu_F(y)\leq &\sum _{j}\int ||\partial_t^\gamma\mathcal{F}\mu_E^j(y,t)||_{L^2(dt)}\,d\mu_F(y)\\
\leq &\sum _{j}\left(\int ||\partial_t^\gamma\mathcal{F}\mu_E^j(y,t)||^2_{L^2(dt)}\,d\mu_F(y)\right)^\frac{1}{2}\\ = &\sum _{j}\left(\sum_{j'}\int ||\partial^\gamma_t\mathcal{F}\mu_E^j(y,t)||^2_{L^2(dt)}\,\mu_F^{j'}(y)\,dy\right)^\frac{1}{2}
\end{aligned}
\end{equation}

To complete the proof, we need the following two lemmas, whose proofs are given at the end of this section.
\begin{lemma}\label{error}
There exists $j_0>0$ such that whenever $j'>j+j_0$,
$$\int ||\partial^\gamma_t\mathcal{F}\mu_E^j(y,t)||^2_{L^2(dt)}\,\mu_F^{j'}(y)\,dy\lesssim_N 2^{-j'N},\ \ \forall\ N>0.$$
\end{lemma}
\begin{lemma}
  \label{mulp}
  $$||\mu_E^j||_{L^p_\alpha}\lesssim_\epsilon 2^{j(\alpha+\frac{d-s_E+\epsilon}{p'})}.  $$
In particular,
$$\int_{|\xi|\lesssim 2^j}|\widehat{\mu_E}(\xi)|^2\,d\xi\lesssim_\epsilon 2^{j(d-s_E+\epsilon)}.$$
\end{lemma}
\vskip.125in
By H\"older's inequality, Theorem \ref{l2} and Lemma \ref{mulp}, for any $j, j'>0$,
\begin{equation}\label{applylocalsmoothing}
\begin{aligned}
\int ||\partial^\gamma_t\mathcal{F}\mu_E^j(y,t)||^2_{L^2(dt)}\,\mu_F^{j'}(y)\,dy\leq &\left( \int ||\partial^\gamma_t\mathcal{F}\mu_E^j(y,t)||^{2q}_{L^2(dt)}\,dy\right)^{\frac{1}{2q}\cdot 2} ||\mu_F^{j'}||_{L^{q'}}\\
\lesssim &||\mu_E^{j}||^2_{L^2_{\alpha_{d,2q,\gamma}-\delta_2(d,2q)}}||\mu_F^{j'}||_{L^{q'}}\\ \lesssim_\epsilon &2^{2j(-\frac{d-1}{2}+\gamma+(d-1)(\frac{1}{2}-\frac{1}{2q})-\delta_2(d,2q)+\frac{d-s_E+\epsilon}{2})}\cdot 2^{j'\frac{d-s_F+\epsilon}{q}}.
\end{aligned}
\end{equation} 
\vskip.125in
Together with Lemma \ref{error}, \eqref{decomp} becomes
\begin{equation}
  \begin{aligned}
    &\int ||\partial_t^\gamma\mathcal{F}\mu(y,t)||_{L^2(dt)}\,d\mu_F(y)\\\leq &\sum _{j}\left(\sum_{j'}\int ||\partial^\gamma_t\mathcal{F}\mu_E^j(y,t)||^2_{L^2(dt)}\,\mu_F^{j'}(y)\,dy\right)^\frac{1}{2}\\\leq & \sum _{j}\left(\sum_{j'<j+j_0}\int ||\partial^\gamma_t\mathcal{F}\mu_E^j(y,t)||^2_{L^2(dt)}\,\mu_F^{j'}(y)\,dy + 2^{-j'N}\right)^\frac{1}{2}\\\lesssim_\epsilon &\sum_j\left(\sum_{j'<j+j_0}2^{2j(-\frac{d-1}{2}+\gamma+(d-1)(\frac{1}{2}-\frac{1}{2q})-\delta_2(d,2q)+\frac{d-s_E+\epsilon}{2})}\cdot 2^{j'\frac{d-s_F+\epsilon}{q}}+2^{-j'N}\right)^\frac{1}{2}\\\lesssim_\epsilon & \sum_{j} 2^{j(\gamma+\frac{1}{2q}-\delta_2(d,2q)+\frac{d}{2}-\frac{s_E-\epsilon}{2}-\frac{s_F-\epsilon}{2q})}.
  \end{aligned}
\end{equation}

Take $q=\frac{(d+1)}{d-1}$ and $\epsilon$ to be small enough, then $\delta_2(d,2q)=\frac{1}{2q}$ and the integral above is finite whenever 
$$s_F>d+1+ (1+2\gamma)\frac{(d+1)}{d-1}-\frac{d+1}{d-1}s_E,  $$
as desired.

\begin{proof}[proof of Lemma \ref{error}]
  By interpolation we may assume $\gamma$ is an integer, then
$$\partial^\gamma_t\mathcal{F}\mu_E^j(y,t) =  \int_{\R^d} \left(\int_{\R} e^{-2\pi i (\Phi(x,y)-t)\tau}\tau^\gamma\,d\tau\right)\mu_E^j(x)\,dx$$
and by Plancherel in $t$,
\begin{equation}\label{x,x',y}
\begin{aligned}
&\int ||\partial^\gamma_t\mathcal{F}\mu_E^j(y,t)||^2_{L^2(dt)}\,\mu_F^{j'}(y)\,dy\\=&\iint \left(\int_{\R}\left(\int_{\R^d} e^{-2\pi i \Phi(x,y)\tau}\tau^\gamma\,\mu_E^j(x)\,dx\right)\,e^{2\pi i t\tau}\,d\tau\right)^2\,dt\,\mu_F^{j'}(y)\,dy\\=&\iint \left(\int_{\R^d} e^{-2\pi i \Phi(x,y)\tau}\tau^\gamma\,\mu_E^j(x)\,dx\right)^2\,d\tau\,\mu_F^{j'}(y)\,dy\\=&\iiiint e^{-2\pi i (\Phi(x,y)-\Phi(x',y))\tau}\tau^{2\gamma}\,\mu_E^j(x)\,\mu_E^j(x')\,\mu_F^{j'}(y)\,dx\,dx'\,dy\,d\tau.
\end{aligned}
\end{equation}
By the definition of Littlewood-Paley decomposition,
$$|\mu^j (x)|= 2^{jd}\left|\int \hat{\phi}(2^j(x-z)\,d\mu(z)\right|\leq C_N 2^{jd} \int (1+2^j|x-z|)^{-N}\,d\mu(z).$$
So it suffices to consider $x,x',y$ in \eqref{x,x',y} in a bounded domain. So it approximately equals
\begin{equation*}\int \left(\int\limits_{|x|,|x'|,|y|\lesssim 1} e^{-2\pi i (\Phi(x,y)\tau-\Phi(x',y)\tau-x\cdot\xi-x'\cdot\xi'-y\cdot\eta)}\,dx\,dx'\,dy\right)\,\widehat{\mu_E^j}(\xi)\,\widehat{\mu_E^j}(\xi')\,\widehat{\mu_F^{j'}}(\eta)\,d\xi\,d\xi'\,d\eta\,\tau^{2\gamma}d\tau.  \end{equation*}

Denote the phase function as $\varphi$, i.e., 
$$\varphi=\Phi(x,y)\tau-\Phi(x',y)\tau-x\cdot\xi-x'\cdot\xi'-y\cdot\eta.$$
Then 
\begin{equation*}
  \begin{aligned}
    \nabla_x\varphi = &\tau\nabla_x\Phi(x,y) -\xi,\\
    \nabla_{x'}\varphi = &-\tau \nabla_{x'}\Phi(x',y)-\xi',\\
    \nabla_y\varphi = &\tau(\nabla_y\Phi(x,y)-\nabla_y\Phi(x',y))-\eta.\\
  \end{aligned}
\end{equation*}
where $|\xi|\approx|\xi'|\approx 2^j$, $|\eta|\approx 2^{j'}$.

Since the Phong-Stein rotational curvature condition \eqref{PS} holds, $|\nabla_x\Phi|\approx|\nabla_{x'}\Phi|\approx|\nabla_y\Phi|\approx 1$.
\vskip.125in
If $|\nabla_y\varphi|\gtrsim 2^{j'}$, the lemma follows by integration by parts. If not, $\tau$ must be $\gtrsim 2^{j'}$. Therefore if $j'>j+j_0$, it follows that $|\nabla_x\varphi|\approx|\nabla_{x'}\varphi|\gtrsim 2^{j'}$ and the Lemma follows by integration by parts.
\end{proof}

\begin{proof}
  [proof of Lemma \ref{mulp}]
  Again by interpolation we may assume $\alpha$ is an even integer. Denote $\phi_j(\cdot)=\phi(\frac{\cdot}{2^j})$. Then $\mu^j=\widehat{\phi_j}*\mu$ and 
$$\Delta^{\frac{\alpha}{2}}\mu_E^j(x)=(\Delta^{\frac{\alpha}{2}}\widehat{\phi_j})*\mu(x).$$
It is easy to see 
$$||\Delta^{\frac{\alpha}{2}}\widehat{\phi_j}||_{L^1}\lesssim  2^{\alpha j}.$$
Therefore $$||\mu_E^j||_{L^1_\alpha}=||(\Delta^{\frac{\alpha}{2}}\widehat{\phi_j})*\mu||_{L^1}\leq ||\Delta^{\frac{\alpha}{2}}\widehat{\phi_j}||_{L^1}||\mu||_{L^1}\lesssim 2^{j\alpha}.$$ 
\vskip.125in
For $L_\alpha^\infty$ norm, notice 
$$|(\Delta^{\frac{\alpha}{2}}\widehat{\phi_j})*\mu(x)|= 2^{j\alpha}\,2^{jd}\left|\int  (\Delta^{\frac{\alpha}{2}}\widehat{\phi})(2^j(x-z))\,d\mu(z)\right|\lesssim_\epsilon 2^{j(d-s_E+\alpha+\epsilon)}  $$
where the last inequality follows by \eqref{Frostman} and the fact $|(\Delta^{\frac{\alpha}{2}}\widehat{\phi})(2^j(x-z))|\lesssim_N (1+2^j|x-z|)^{-N}$. This upper bound is uniformly in $x$ so we have 
$$||\mu_E^j||_{L^\infty_\alpha}\lesssim_\epsilon 2^{j(d-s_E+\alpha+\epsilon)}$$
and the lemma follows by interpolation.
\end{proof}

\section{Proof of Theorem \ref{slicingthm}}
 In the last section we proved that $\nu_y(t)$ on 
 $$\Delta^y_{\Phi}(E)=\{\Phi(x,y):x\in E\},$$
 defined in \eqref{pinmeasure}, is in fact an $L^2$ function for $\mu_F$-a.e. $y\in F$. With this in mind one can define the slicing measure $\mu_E^{y,t}$ on $\{x\in E:\Phi(x,y)=t\}$ as
$$\int f\,d\mu_E^{y,t}= \iint e^{-2\pi i (\Phi(x,y)-t)\tau}\,d\tau\,f(x)d\mu_E(x).  $$
Notice $\mu_E^{y,t}$ is a well-defined measure for $\mu_F\times\mathcal{L}^1$-a.e. $(y,t)\in F\times\R$ because by definition
$$\int d\mu_E^{y,t} = \nu_y(t)<\infty  $$
for $\mu_F\times\mathcal{L}^1$-a.e. $(y,t)\in F\times\R$. Also because
$$\int\nu_y(t)\,dt=\int\,d\mu_E(x)=1,$$
$\{t\in\R:\int \mu_E^{y,t}>0\}$ has positive Lebesgue measure for $\mu_F$-a.e. $y\in F$.
\vskip.125in
Define the $\sigma$-energy of $\mu_E^{y,t}$ as
$$I_{\sigma}(\mu_E^{y,t})=\iint |u-v|^{-\sigma}\,d\mu_E^{y,t}(u)\,d\mu_E^{y,t}(v)=c_{\sigma,d}\int |\widehat{\mu_E^{y,t}}(\xi)|^2|\xi|^{-d+\sigma}\,d\xi.  $$
It is known that the finiteness of $I_{\sigma}(\mu_E^{y,t})$ implies the Hausdorff dimension of the support of $\mu_E^{y,t}$ is at least $\sigma$ (see, e.g. \cite{Mat95}). Thus to prove Theorem \ref{slicingthm}, it suffices to show
\begin{equation}
  \label{goalofslicing}
  \iint I_{\sigma}(\mu_E^{y,t})\,dt\,d\mu_F(y)<\infty
\end{equation}
whenever $\sigma<s_E+\frac{d+1}{d-1}s_F-d$.

\vskip.125in
By the definition of $I_{\sigma}(\mu_E^{y,t})$,
\begin{equation}\label{1}
\begin{aligned}
&\iint I_{\sigma}(\mu_E^{y,t})\,dt\,d\mu_F(y)\\
=&  \iiint|\widehat{\mu_E^{y,t}}(\xi)|^2 |\xi|^{-d+\sigma}\,d\xi\,dt\,d\mu_F(y)\\=&\iiint\left|\int \int e^{-2\pi i (\Phi(x,y)-t)\tau}\,d\tau\,e^{-2\pi i x\cdot\xi}d\mu_E(x)\right|^2 |\xi|^{-d+\sigma}\,d\xi\,dt\,d\mu_F(y).
\end{aligned} 
\end{equation}
Denote $\mu_E^\xi = e^{-2\pi i x\cdot\xi}\mu_E$. By Plancherel in $t$, \eqref{1} becomes
\begin{equation*}
\begin{aligned}
&\iiint\left|\int e^{-2\pi i \Phi(x,y)\tau}\,d\mu^\xi_E(x)\right|^2 |\xi|^{-d+\sigma}\,d\xi\,d\tau\,d\mu_F(y)\\=&\int\cdots\int e^{-2\pi i ((\Phi(x,y)-\Phi(x',y))\tau}\,d\mu^\xi_E(x)\,d\mu^\xi_E(x')\,|\xi|^{-d+\sigma}d\xi\,d\tau\,d\mu_F(y)\\=&\int\left(\int e^{-2\pi i ((\Phi(x,y)-\Phi(x',y))\tau-x\cdot\eta-x'\cdot\eta'-y\cdot\zeta)}\,dx\,dx'\,dy\right)\widehat{\mu^\xi_E}(\eta)\,\widehat{\mu^\xi_E}(\eta')\widehat{\mu_F}(\zeta)\,|\xi|^{-d+\sigma}\,d\xi\,d\eta\,d\eta'\,d\zeta\,d\tau.
\end{aligned}
\end{equation*}
Since $E, F$ are both compact sets, $x,x',y$ above lie in a bounded domain. Denote $\varphi$ as the phase function in the integral above, i.e., 
$$\varphi=(\Phi(x,y)-\Phi(x',y))\tau-x\cdot\eta-x'\cdot\eta'-y\cdot\zeta.  $$
Then
\begin{equation*}
  \begin{aligned}
    \nabla_x\varphi = &\tau\nabla_x\Phi(x,y) - \eta,\\
    \nabla_{x'}\varphi = &-\tau \nabla_{x'}\Phi(x',y)-\eta',\\
    \nabla_y\varphi = &\tau(\nabla_y\Phi(x,y)-\nabla_y\Phi(x',y))-\zeta.\\
  \end{aligned}
\end{equation*}
Since the Phong-Stein rotational curvature condition \eqref{PS} holds, $|\nabla_x\Phi|\approx|\nabla_{x'}\Phi|\approx|\nabla_y\Phi|\approx 1$. 
\vskip.125in
When $|\nabla\varphi|=0$, $\tau\approx |\eta|\approx|\eta'|\gtrsim|\zeta|$. Thus it suffices to consider the domains where $|\eta|\approx|\eta'|\approx2^{j_1}$, $|\zeta|\approx 2^{j_2}$, $j_1\geq j_2$.
\vskip.125in
For each fixed $\xi$, denote
$$\mu^{\xi, j_1}_E =\left(\widehat{\mu^\xi_E}(\cdot)\phi(\frac{\cdot}{2^{j_1}})\right)^\vee= \left(\widehat{\mu_E}(\cdot+\xi)\phi(\frac{\cdot}{2^{j_1}})\right)^\vee.  $$

The discussion above implies that, to show the finiteness of \eqref{1}, it suffices to consider
\begin{equation*}
\begin{aligned}
&\sum_{j_1\geq j_2}\int\int \left(\int\left|\int \int e^{-2\pi i (\Phi(x,y)-t)\tau}\,d\tau\,\mu^{\xi,j_1}_E(x)\,dx\right|^2 \,dt\right)\,\mu^{j_2}_F(y)\,dy\, |\xi|^{-d+\sigma}\,d\xi \\=& \sum_{j_1\geq j_2}\int\left(\int \left(\int |\mathcal{F}\mu^{\xi,j_1}_E(y,t)|^2 \,dt\right)\,\mu^{j_2}_F(y)\,dy\right) \,|\xi|^{-d+\sigma}\,d\xi.
\end{aligned}
\end{equation*}
Fix $\xi$, apply H\"older's inequality in $dy$ with $q>1$. It is less than or equal to
\begin{equation}\label{2}
 \sum_{j_1\geq j_2}||\mu^{j_2}_F||_{L^{q'}}\,\int \left(\int\left(\int |\mathcal{F}\mu^{\xi,j_1}_E(y,t)|^2 \,dt\right)^{\frac{1}{2}\cdot 2q}\,dy\right)^{\frac{1}{2q}\cdot2}\,|\xi|^{-d+\sigma}\,d\xi. 
\end{equation}
\vskip.125in
Then by Lemma \ref{mulp} and the $L^2\rightarrow L^p$ sharp local smoothing estimate, with $2q=\frac{2(d+1)}{d-1}$, this integral is bounded above by
\begin{equation}\label{4}
\begin{aligned}
&\sum_{j_1\geq j_2}2^{j_2\frac{d-s_F+\epsilon}{q}}\int||\mu^{\xi,j_1}_E(y,t)||^2_{L^2_{\alpha_{d,2q,0}-\frac{1}{2q}}}\,|\xi|^{-d+\sigma}\,d\xi\\=&\sum_{j_1}2^{j_1\frac{d-s_F+\epsilon}{q}}\int||\mu^{\xi,j_1}_E(y,t)||^2_{L^2_{\alpha_{d,2q,0}-\frac{1}{2q}}}\,|\xi|^{-d+\sigma}\,d\xi.
\end{aligned}
\end{equation}
\vskip.125in
So to prove \eqref{goalofslicing}, it suffices to show \eqref{4} is finite.
\vskip.125in
By the definition of $\mu^{\xi,j_1}_E$, 
$$\widehat{\mu^{\xi,j_1}_E}(\eta) =\widehat{\mu_E}(\eta+\xi)\phi(\frac{\eta}{2^{j_1}}).$$
So if we take $2q=\frac{2(d+1)}{d-1}$,
\begin{equation}\label{5}
\begin{aligned}
\int||\mu^{\xi,j_1}_E(y,t)||^2_{L^2_{\alpha_{d,2q}-\frac{1}{2q}}}\,|\xi|^{-d+\sigma}\,d\xi = &\int\int_{|\eta|\approx{2^{j_1}}} |\widehat{\mu^\xi_E}(\eta)|^2\,|\eta|^{-\frac{d(d-1)}{d+1}}d\eta\,|\xi|^{-d+\sigma}\,d\xi\\\approx&\,2^{-j_1\frac{d(d-1)}{d+1}}\int\int_{|\eta|\approx{2^{j_1}}} |\widehat{\mu_E}(\eta+\xi)|^2\,d\eta\,|\xi|^{-d+\sigma}\,d\xi\\=&\,2^{-j_1\frac{d(d-1)}{d+1}}\left(\int_{|\xi|\lesssim{2^{j_1}}}\int_{|\eta|\approx{2^{j_1}}}+\sum_{j_3\geq j_1}\int_{|\xi|\approx{2^{j_3}}}\int_{|\eta|\approx{2^{j_1}}}\right)\\=&\,2^{-j_1\frac{d(d-1)}{d+1}}(I+II).
\end{aligned}
\end{equation}
\vskip.125in
For $I$, since $|\xi|\lesssim{2^{j_1}}$, $|\eta|\approx{2^{j_1}}$, we have $|\eta+\xi|\lesssim{2^{j_1}}\approx|\eta|$. Then by changing variables $\eta'=\eta+\xi$ and Lemma \ref{mulp},
\begin{equation*}
\begin{aligned}
I&=\int_{|\xi|\lesssim{2^{j_1}}}\int_{|\eta|\approx{2^{j_1}}}|\widehat{\mu_E}(\eta+\xi)|^2\,d\eta\,|\xi|^{-d+\sigma}\,d\xi\\&\lesssim\int_{|\xi|\lesssim{2^{j_1}}}\int_{|\eta'|\lesssim{2^{j_1}}}|\widehat{\mu_E}(\eta')|^2\,d\eta'\,|\xi|^{-d+\sigma}\,d\xi\\&\lesssim_\epsilon 2^{j_1(d-s_E+\epsilon+\sigma)}.
\end{aligned}
\end{equation*}
\vskip.125in
For $II$, since $|\xi|\approx 2^{j_3}$, $|\eta|\approx 2^{j_1}$, $j_3\geq j_1$, it follows $|\eta+\xi|\lesssim{2^{j_3}}\approx|\xi|$. Then by changing variables $\xi'=\eta+\xi$ and Lemma \ref{mulp},

\begin{equation*}
\begin{aligned}II&=\int_{|\xi|\approx 2^{j_3}}\int_{|\eta|\approx 2^{j_1}}|\widehat{\mu_E}(\eta+\xi)|^2\,d\eta\,|\xi|^{-d+\sigma}\,d\xi\\&\lesssim  2^{j_3(-d+\sigma)}\int_{|\eta|\approx 2^{j_1}}\left(\int_{|\xi|\approx 2^{j_3}}|\widehat{\mu_E}(\eta+\xi)|^2\,d\xi\right)\,d\eta\\&\lesssim  2^{j_3(-d+\sigma)}\int_{|\eta|\approx 2^{j_1}}\left(\int_{|\xi'|\approx 2^{j_3}}|\widehat{\mu_E}(\xi')|^2\,d\xi'\right)\,d\eta\\&\lesssim_\epsilon 2^{j_3(-d+\sigma)}\,2^{j_3(d-s_E+\epsilon)}\,2^{j_1d}\\&=\,2^{j_3(\sigma-s_E+\epsilon)}\,2^{j_1d}.
\end{aligned}
\end{equation*}

\vskip.125in

These estimates of $I, II$ give us an upper bound of \eqref{5}. Thus \eqref{4} is bounded above, with $2q=\frac{2(d+1)}{d-1}$, by
$$\sum_{j_1}2^{j_1(d-s_F+\epsilon)\frac{(d-1)}{d+1}}\,2^{-j_1\frac{d(d-1)}{d+1}}\left(2^{j_1(d-s_E+\epsilon+\sigma)}+\sum_{j_1\leq j_3}\,2^{j_3(\sigma-s_E+\epsilon)}\,2^{j_1d}\right).  $$

Since $\sigma$ is the expected dimension of slicing measures of $E$, we may assume $\sigma<s_E$, then the sum is no bigger than
$$\sum_{j_1}2^{j_1(d-s_F+\epsilon)\frac{(d-1)}{d+1}}\,2^{-j_1\frac{d(d-1)}{d+1}}\,2^{j_1(d-s_E+\epsilon+\sigma)},  $$
which is finite whenever $\sigma<s_E+\frac{d-1}{d+1}s_F-d$ and $\epsilon$ is sufficiently small.

\vskip.25in
\section{Proof of Theorem \ref{sharpl2} and Theorem \ref{sharplp}}
We prove Theorem \ref{sharpl2} first. If we run the proof of Theorem \ref{main}, without taking a specific value of $q$ in the last step, it follows that 
\begin{equation}\label{Minkowski}\dH( \{y\in V: |\Delta_{\Phi}^y(E)| = 0 \})\leq qd+1-2q\delta_2(d,2q)-q\dH(E),\ \forall q>1.
\end{equation}
\vskip.125in

For $(i)$, if $\delta_2(d,2q)>\frac{1}{2q}$, the right hand side of \eqref{Minkowski} is negative when $\dH(E)=d$, which means
$$ \{y\in V: |\Delta_{\Phi}^y(E)| = 0 \}=\emptyset.$$

However, one can take $E=\{x:\Phi(x,y_0)\in A\}$, where $y_0\in V$ is fixed and $A\subset\R$ has Hausdorff dimension $1$ but Lebesgue measure $0$. Then $\dH(E)=d$, while $|\Delta_{\Phi}^{y_0}(E)| = |A|=0$. Contradiction.
\vskip.125in
For $(ii)$, say $d=2k+1$, one can take $\Phi(x,y)$ to be the $(k+1,k)$ - Minkowski distance between $x$ and $y$, i.e.
$$\Phi(x,y) = ||x-y||_{k+1,k}=\sum_{i=1}^{d}(-1)^{i+1}(x_i-y_i)^2$$
and 
$$E= A\times \{(t,t)\in\R^2:t\in\R \}^k,$$ where $A\subset\R$ has Hausdorff dimension $1$ but Lebesgue measure $0$. Then for any point $$y\in\{(y_1,\dots,y_{2k+1})\in\R^{2k+1}:y_{2j}=y_{2j+1}, 1\leq j\leq k \},$$whose Hausdorff dimension is $\frac{d+1}{2}$, it is easy to see that $|\Delta_{\Phi}^y(E)|=|y_1|\cdot|A| = 0$. This means 
$$\dH( \{y\in\Omega: |\Delta_{\Phi}^y(E)| = 0 \})\geq \frac{d+1}{2}.$$
Compared with \eqref{Minkowski}, $\delta_2(d,2q)$ has to be less than or equal to $\frac{1}{2}(d-1)(\frac{1}{2}-\frac{1}{2q})$.
\vskip.125in
For $(iii)$, say $d=2k$, one can take 
$$\Phi(x,y)= x_1y_1+\sum_{i=2}^{d}(-1)^{i}(x_i-y_i)^2$$
and $$E=A\times [1,2]\times \{(t,t)\in\R^2:t\in\R \}^{k-1},$$
where $A\subset\R$ has Hausdorff dimension $1$ but Lebesgue measure $0$. Then for any 
$$y\in\{(y_1,\dots,y_{2k})\in\R^{2k}: y_{2j-1}=j_{2j}, 2\leq j\leq k\},$$whose Hausdorff dimension is $\frac{d}{2}$, one can check that $|\Delta_{\Phi}^y(E)| = 0$. This means
$$\dH( \{y\in\Omega: |\Delta_{\Phi}^y(E)| = 0 \})\geq \frac{d}{2}.$$
Compare with \eqref{Minkowski}, $\delta_2(d,2q)$ has to be less than or equal to $(d-1)(\frac{1}{2}-\frac{1}{2q})$.
\vskip.25in
For Theorem \ref{sharplp}, we apply Conjecture \ref{Sog} instead of Theorem \ref{l2} in the proof of Theorem \ref{main}. More precisely, we apply H\"older's inequality
\begin{equation*}
\begin{aligned}\left|\int ||\partial^\gamma_t\mathcal{F}\mu_E^j(y,t)||^2_{L^2(dt)}\,\mu_F^{j'}(y)\,dy\right|\leq &\int ||\partial^\gamma_t\mathcal{F}\mu_E^j(y,t)||^2_{L^p(dt)}\,|\mu_F^{j'}(y)|\,dy\\\leq&||\partial^\gamma_t\mathcal{F}\mu_E^j(y,t)||^2_{L^{p}(V\times[1,2])}\,||\mu_F^{j'}||_{L^{(p/2)'}(V)}\\\lesssim &||\mu_E^{j}||^2_{L^{p}_{\alpha_{d,p,\gamma}-\delta_{p}(d,p)}}||\mu_F^{j'}||_{L^{(p/2)'}}\\ \lesssim&_\epsilon 2^{2j(-\frac{d-1}{2}+\gamma+(d-1)(\frac{1}{2}-\frac{1}{p})-\delta_{p}(d,p)+\frac{d-s_E+\epsilon}{p'})}\cdot 2^{j'\frac{2(d-s_F+\epsilon)}{p}}.
\end{aligned}
\end{equation*}

It follows that, for any $p>2$,
$$\dH( \{y\in V: |\Delta_{\Phi}^y(E)| = 0 \})\leq \frac{p}{p'}d+1-p\delta_p(d,p)-\frac{p}{p'}\dH(E).$$

Then Theorem \ref{sharplp} follows with the same counterexamples in the proof of Theorem \ref{sharpl2}.

\bibliographystyle{abbrv}
\bibliography{/Users/MacPro/Dropbox/Academic/paper/mybibtex.bib}

\begin{thebibliography}{10}

\bibitem{BD15}
J.~Bourgain and C.~Demeter.
\newblock The proof of the {$l^2$} decoupling conjecture.
\newblock {\em Ann. of Math. (2)}, 182(1):351--389, 2015.

\bibitem{Erd05}
M.~B. Erdogan.
\newblock A bilinear {F}ourier extension theorem and applications to the
  distance set problem.
\newblock {\em Int. Math. Res. Not.}, (23):1411--1425, 2005.

\bibitem{EIT11}
S.~Eswarathasan, A.~Iosevich, and K.~Taylor.
\newblock Fourier integral operators, fractal sets, and the regular value
  theorem.
\newblock {\em Adv. Math.}, 228(4):2385--2402, 2011.

\bibitem{Fal85}
K.~J. Falconer.
\newblock On the {H}ausdorff dimensions of distance sets.
\newblock {\em Mathematika}, 32(2):206--212, 1985.

\bibitem{Hor71}
L.~H{\"o}rmander.
\newblock Fourier integral operators. {I}.
\newblock {\em Acta Math.}, 127(1-2):79--183, 1971.

\bibitem{IL16}
A.~Iosevich and B.~Liu.
\newblock Falconer distance problem, additive energy and cartesian products.
\newblock {\em Ann. Acad. Sci. Fenn. Math.}, 41(2):579--585, 2016.

\bibitem{ITU16}
A.~Iosevich, K.~Taylor, and I.~Uriarte-Tuero.
\newblock Pinned geometric configurations in euclidean space and riemannian
  manifolds.
\newblock {\em https://arxiv.org/pdf/1610.00349v1.pdf}, 2016.

\bibitem{Mat95}
P.~Mattila.
\newblock {\em Geometry of Sets and Measures in Euclidean Spaces: Fractals and
  Rectifiability}, volume~44 of {\em Cambridge Studies in Advanced
  Mathematics}.
\newblock Cambridge University Press, Cambridge, 1995.

\bibitem{Mat14}
P.~Mattila.
\newblock Recent progress on dimensions of projections.
\newblock In {\em Geometry and analysis of fractals}, volume~88 of {\em
  Springer Proc. Math. Stat.}, pages 283--301. Springer, Heidelberg, 2014.

\bibitem{Mat15}
P.~Mattila.
\newblock {\em Fourier analysis and Hausdorff dimension}, volume 150.
\newblock Cambridge University Press, 2015.

\bibitem{MS97}
W.~P. Minicozzi, II and C.~D. Sogge.
\newblock Negative results for {N}ikodym maximal functions and related
  oscillatory integrals in curved space.
\newblock {\em Math. Res. Lett.}, 4(2-3):221--237, 1997.

\bibitem{MSS93}
G.~Mockenhaupt, A.~Seeger, and C.~D. Sogge.
\newblock Local smoothing of {F}ourier integral operators and
  {C}arleson-{S}j\"olin estimates.
\newblock {\em J. Amer. Math. Soc.}, 6(1):65--130, 1993.

\bibitem{Orp14}
T.~Orponen.
\newblock Slicing sets and measures, and the dimension of exceptional
  parameters.
\newblock {\em J. Geom. Anal.}, 24(1):47--80, 2014.

\bibitem{Orp17}
T.~Orponen.
\newblock On the distance sets of {A}hlfors-{D}avid regular sets.
\newblock {\em Adv. Math.}, 307:1029--1045, 2017.

\bibitem{PS00}
Y.~Peres and W.~Schlag.
\newblock Smoothness of projections, {B}ernoulli convolutions, and the
  dimension of exceptions.
\newblock {\em Duke Math. J.}, 102(2):193--251, 2000.

\bibitem{Shm16}
P.~Shmerkin.
\newblock On distance sets, box-counting and ahlfors-regular sets.
\newblock {\em arXiv:1604.00308}, 2016.

\bibitem{shm17}
P.~Shmerkin.
\newblock On the hausdorff dimension of pinned distance sets.
\newblock {\em arXiv preprint arXiv:1706.00131}, 2017.

\bibitem{Sog91}
C.~D. Sogge.
\newblock Propagation of singularities and maximal functions in the plane.
\newblock {\em Invent. Math.}, 104(2):349--376, 1991.

\bibitem{Sog93}
C.~D. Sogge.
\newblock {\em Fourier integrals in classical analysis}, volume 105 of {\em
  Cambridge Tracts in Mathematics}.
\newblock Cambridge University Press, Cambridge, 1993.

\bibitem{Wol99}
T.~Wolff.
\newblock Decay of circular means of {F}ourier transforms of measures.
\newblock {\em Internat. Math. Res. Notices}, (10):547--567, 1999.

\bibitem{Wol00}
T.~Wolff.
\newblock Local smoothing type estimates on {$L^p$} for large {$p$}.
\newblock {\em Geom. Funct. Anal.}, 10(5):1237--1288, 2000.

\end{thebibliography}

\end{document}